\documentclass[graybox]{svmult}

\usepackage{mathptmx}       
\usepackage{helvet}         
\usepackage{courier}        
\usepackage{type1cm}        
\usepackage{makeidx}         
\usepackage{graphicx}        
\usepackage{multicol}        

\usepackage{amssymb,amsfonts,amscd,amsmath}

\input amssym.def

\input amssym.tex

\makeindex             

\newtheorem{thm}{Theorem}[section]

\newtheorem{rek}[thm]{Remark}

\newcommand\bi{\begin{itemize}}

\newcommand\ei{\end{itemize}}

\newcommand\ben{\begin{enumerate}}

\newcommand\een{\end{enumerate}}

\newcommand{\gep}{\epsilon}  

\newcommand\be{\begin{equation}}

\newcommand\ee{\end{equation}}

\newcommand\bea{\begin{eqnarray}}

\newcommand\eea{\end{eqnarray}}

\numberwithin{subsubsection}{subsection}

\begin{document}

\title*{Finding and Counting MSTD sets}

\author{Geoffrey Iyer, Oleg Lazarev, Steven J. Miller and Liyang Zhang}
\institute{Department of Mathematics, University of Michigan, geoff.iyer@gmail.com; Department of Mathematics, Princeton University, olazarev@Princeton.EDU; Department of Mathematics and Statistics, Williams College, sjm1@williams.edu (Steven.Miller.MC.96@aya.yale.edu) and lz1@williams.edu.}

%
%
\maketitle

\abstract{We review the basic theory of More Sums Than Differences (MSTD) sets, specifically their existence, simple constructions of infinite families, the proof that a positive percentage of sets under the uniform binomial model are MSTD but not if the probability that each element is chosen tends to zero, and `explicit' constructions of large families of MSTD sets. We conclude with some new constructions and results of generalized MSTD sets, including among other items results on a positive percentage of sets having a given linear combination greater than another linear combination, and a proof that a positive percentage of sets are $k$-generational sum-dominant (meaning $A$, $A+A$, $\dots$, $kA = A + \cdots +A$ are each sum-dominant).
}



\section{Introduction}

Many of the most important questions in additive number theory can be cast as questions about sums or differences of sets, where the sumset of $A$ and $B$ is  \be A+B \ = \ \{a+b: A \in A, b \in B\}\ee and the difference set is \be A-B \ = \ \{a-b: a\in A, b \in B\}. \ee To see this, let $\mathcal{P}$ be the set of primes and $\mathcal{N}_k$ (respectively $\mathcal{N}_k'$) be the set of $k$\textsuperscript{th} powers of integers (respectively non-negative integers).


\ben

\item The famous Goldbach problem is to prove that every even number may be written as the sum of two primes; we may interpret this as saying that the even numbers are contained in $\mathcal{P}+\mathcal{P}$. While this is still open, we do know that all sufficiently large odd numbers are the sum of three primes. While sufficiently large means greater than $10^{1000}$ here,  we may remove `sufficiently large' if we assume the Generalized Riemann Hypothesis \cite{DETZ}.

\item Another example is Waring's problem, which says for each integer $k$ there is an integer $s$ such that every positive integer is a sum of at most $s$ perfect $k$\textsuperscript{th} powers. In other words, there is an $s$ (depending on $k$) such that $\mathcal{N}_k + \cdots + \mathcal{N}_k$ (where there are $s$ sums) contains all positive integers. While the optimal $s$ for a given $k$ is not known, it is known that for each $k$ there does exist a finite $s$ (see for instance \cite{Na1}).

\item Fermat's Last Theorem (proved in \cite{Wi,TW}) states that if $n\ge 3$ and $x,y,z$ are integers, then the only solutions to $x^n+y^n=z^n$ have $xyz=0$. After some simple algebra we see it suffices to consider the case when $x,y$ and $z$ are all positive, and Fermat's Last Theorem is just the statement that $(\mathcal{N}_n' + \mathcal{N}_n') \cap \mathcal{N}_n'$ is empty for $n \ge 3$.

\een

The three examples above all involve determining what elements are in sums of sets; it is also interesting to see how often a given element is represented in a sum. For example, the Twin Prime Conjecture is the assertion that there are infinitely many primes differing by 2; this is equivalent to how often 2 is obtained in $\mathcal{P}_x - \mathcal{P}_x$, where $\mathcal{P}_x$ is the truncated set of primes at most $x$.

As the topic of sum sets and difference sets is so vast, in this survey article we restrict ourselves to an interesting class of questions where there has been significant progress in recently years. Given a finite set of integers $A$, we may look at $A+A$ and $A-A$. The most natural question to ask is: As we vary $A$ over a family of sets, how often is the cardinality of $A+A$ larger than $A-A$? Denoting the size of a set $S$ by $|S|$, for $|A|$ large we expect a typical $A$ to have $|A+A| < |A-A|$. This is because while the diagonal pairs $(a,a)$ contribute a new sum to $A+A$ for each $a$ but only one difference (namely 0) to $A-A$, addition is commutative while subtraction is not. This means that for the larger collection of pairs of distinct elements $(a,a')$ we have $a+a' = a'+a$ but $a-a' \neq a'-a$. We see a typical pair contributes two differences to $A-A$ but only one sum to $A+A$. Using such logic, one expects sets with $|A+A| > |A-A|$ to be rare.

If $|A+A| > |A-A|$, we say $A$ is a \emph{sum-dominated} set or a \emph{More Sums Than Differences (MSTD)} set, while if $|A+A| = |A-A|$ we say $A$ is \emph{balanced}, and if $|A+A| < |A-A|$ then $A$ is \emph{difference-dominated}.  The purpose of this article is to describe results in the following areas.

\begin{enumerate}

\item \emph{Non-probabilistic constructions of MSTD sets.} In this section we summarize some of the early constructions of MSTD sets, paying special attention to the limitation of these techniques in determining whether or not a \emph{typical} set is sum-dominated.

\item \emph{A positive percentage of sets are MSTD sets.} Here we discuss the papers of Martin and O'Bryant \cite{MO} and Zhao \cite{Zh2}, which show that a very small, but positive, percentage of all sets are sum-dominated.

\item \emph{When a `typical' subset is difference-dominated.} If we choose our subsets of $\{0,\dots,n-1\}$ from the uniform model, so that each of the $2^n$ possible subsets is equally likely to be chosen, then the previous section shows a positive percentage of subsets are sum-dominated. The situation is drastically different if we sample differently. We describe the results of Hegarty and Miller \cite{HM1}, who showed that if each element from $\{0,\dots,n-1\}$ is chosen with probability $p(n)$ and $\lim_{n\to\infty} p(n) = 0$, then in the limit almost all subsets are difference-dominated.

\item \emph{Explicit constructions of large families of MSTD sets.} The methods of \cite{MO,Zh2} are probabilistic, and do not yield explicit families of MSTD sets. Miller, Orosz and Scheinerman \cite{MOS} gave an explicit construction of a large family of subsets of $\{0,\dots,n-1\}$ that are MSTD sets, specifically one whose cardinality is at least $C/n^4$ for some $C>0$; later Zhao \cite{Zh1} gave a different construction yielding $C'/n$ with $C'>0$. We describe these constructions and generalizations; for example, Miller, Pegado and Robinson \cite{MPR} show that the density of sets $A \subset \{0,\dots,n-1\}$ with $|A+A+A+A| > |A+A-A-A|$ is at least $C''/n^r$, where $r = \frac16 \log_2(256/255) \le .001$.

\item \emph{Generalized MSTD Sets.} A set $A$ is a $k$-generational sum-dominant set if $A$, $A+A$, $\dots$, $kA = A + \cdots +A$ are each sum-dominant. Iyer, Lazarev, Miller and Zhang \cite{ILMZ} proved that a positive percentage of sets are $k$-generational for any positive $k$, but no set is $k$-generational for all $k$. Their construction uses a result of interest in its own right, namely that if we are given any legitimate order of linear combinations of sums and differences of $A$ of the same length\footnote{Note that $A+A+A-A = -(A-A-A-A)$; thus we might as well assume any linear combination has at least as many sums of $A$ as differences of $A$.}, a positive percentage of $A$ have the cardinalities of these combinations in the desired ordering. Such a result was expected from the work of Miller, Orosz and Scheinerman \cite{MOS}, who showed if there exists one set satisfying the ordering then there exists a large, explicitly constructible family of sets satisfying the condition. In \cite{ILMZ} the needed set for the induction is found, and instead of appealing to results from \cite{MOS}, the authors modify the arguments of \cite{MO} in order to obtain a positive percentage.

\end{enumerate}

The above list of topics is not meant to be definitive or exhaustive, but rather to highlight some of the many results in the field. There are numerous generalizations to other linear combinations of sets, as well as related problems in Abelian groups, that can be handled with these methods. We strongly urge the reader to consult the references for full details and statements of related, open questions.

\ \\

Miller thanks Mel Nathanson who, through books and conversations, helped introduce him to this exciting subject, his collaborators Peter Hegarty, Brooke Orosz, Sean Pegado, Luc Robinson and Dan Scheinerman for the insights gleaned from our studies, and the participants of various CANT Conferences (especially Greg Martin, Kevin O'Bryant and Jonathan Sondow) for many enlightening conversations; all authors thank the participants of SMALL 2011 for helpful conversations and discussions. The first, second and fourth named authors were supported by NSF grants DMS0850577 and Williams College; the third named author was partially supported by NSF grant DMS0970067.


\section{Non-probabilistic Constructions of MSTD sets.}

In \cite{Na2}, Nathanson wrote \emph{``Even though there
exist sets $A$ that have more sums than differences, such sets should be rare, and it must be true with the right way of counting that the vast majority of sets satisfies $|A-A| > |A+A|$.''} Support for this view can be found in the length of the search required to find the first MSTD set. Conway is said to have found $\{0, 2, 3, 4, 7, 11, 12, 14\}$ in the 1960s, while Marica \cite{Ma} in 1969 gave $\{0, 1, 2, 4, 7, 8, 12, 14, 15\}$ and Freiman and Pigarev \cite{FP} found $\{0, 1, 2, 4, 5$, $9, 12, 13$, $14, 16, 17$, $21, 24, 25, 26, 28, 29\}$ in 1973. See also the papers by Ruzsa \cite{Ru1,Ru2,Ru3}.

How hard is it to find such sets? A simple calculation shows that if $B = \alpha A + \beta$, then $|A+A| = |B+B|$ and $|A-A| = |B-B|$; thus we might as well assume $0$ is in our subset. The number of subsets of $\{0,\dots,14\}$ that include $0$ is $2^{14} = 16,384$. This is easily searchable by computer, though a little out of the range of even the most patient of mathematicians; the only MSTD set found is the one already mentioned. Even Freiman and Pigarev's example can be found by a brute force within a reasonable time, as $2^{29} = 536,870,912$.

While there are many constructions of MSTD sets, most of these constructions give a vanishingly small percentage of sets to be sum-dominated. Specifically, while there are $2^{n+1}$ subsets of $\{0,1,\dots,n\}$, these methods often give only on the order of $2^{n/2}$ (or worse) subsets that are MSTD.

For example, one way to generate an infinite family of MSTD sets from one known MSTD set is through the \emph{base expansion method}. Let $A$ be an MSTD set, and let $A_{k;m} = \{\sum_{i=1}^k a_i m^{i-1}: a_i \in A\}$. If $m$ is sufficiently large, then $|A_{k;m} \pm A_{k;m}| = |A\pm A|^k$. We thus obtain an infinite family of MSTD sets, and, so long as $|A+A| > 1$, we can have arbitrarily many more sums than differences. Unfortunately, as $m$ is large, the percentage of subsets created that are sum-dominated is exponentially small. We thus discuss other constructions (though this method will play an important role in proving many of the theorems in \S\ref{sec:generalizedMSTDsets}).

It is very easy to create balanced sets, and many constructions of MSTD sets take advantage of this. First, note that if $A$ is an arithmetic progression then $A$ is balanced. To see this, letting $A = \{0,1,\dots,n\}$ we find $A+A=\{0,1,\dots, 2n\}$ and $A-A=\{-n,\dots,n\}$ so $|A+A| = |A-A| = 2n+1$. Another way to create a balanced set is to take a set symmetric with respect to a number (which need not be in the set); this means that there is a number $a^\ast$ such that $A = a^\ast - A$ (this implies $A+A = a^\ast + A - A$, so $|A+A| = |A-A|$). Note arithmetic progressions are a special case, with $a^\ast = n/2$. Nathanson \cite{Na3} gives constructions of MSTD sets using this idea. He creates infinite families by adjoining one number to a symmetric set which is a small permutation of a generalized arithmetic progression. Numerous examples and explicit constructions are given in \cite{Na3}; we state the first.

\begin{thm}[Nathanson \cite{Na3}] Let $m,d$, and $k$ be integers with $m \ge 4$, $1 \le d \le m-1$, $d \neq m/2$, and $k \ge 3$ if $d < m/2$ and $k\ge 4$ if $d > m/2$. Let $B = \{0,1,\dots,m-1\} \setminus \{d\}$, $L = \{m-d, 2m-d, \dots, km-d\}$, $a^\ast = (k+1)m-2d$, and $A^\ast = B \cup L \cup (a^\ast - B)$. Then $A = A^\ast \cup \{m\}$ is an MSTD set. \end{thm} 

How large of a family is this? We have three parameters at our disposal: $m, d$ and $k$. Note $A \subset \{0,\dots, (k+1)m-2d\}$. Given some $n$, look at all triples $(m,d,k)$ such that $(k+1)m-2d \le n$; this will be an upper bound for the number of MSTD sets generated by the theorem that live in $\{0,1,\dots,n\}$ (it will be the actual number if we show all the sets are distinct). As we also need $k$ to be at least three, we obtain an upper bound by counting all pairs $(k,m)$ with $km \le n$ (which is trivially at most $n^2$) and noting that we have $m \le n$ choices of $d$ for each pair.  Thus this method generates at most $n^3$ subsets of $\{0,1,\dots,n\}$ being MSTD sets, which is a vanishingly small fraction in the limit. The paucity of this family is due to how explicit the construction is -- everything is completely deterministic and at each stage there is only one option.

We conclude our discussion on constructions of MSTD sets and families of MSTD sets with a result of Hegarty \cite{He}. He proved

\begin{thm}[Hegarty \cite{He}] There are no MSTD subsets of the integers of size seven. Up to linear transformations the only set of size 8 is $\{0,2,3,4,7,11,12,14\}$. \end{thm}  

We paraphrase (slightly) from \cite{He} the description of the proof. Let $A = \{a_n = 0, a_{n-1},\dots,a_1\}$, and represent the $n-1$ differences $a_i - a_{i+1}$ as $\overrightarrow{e}_i$ (the $i$\textsuperscript{th} standard basis vector in $\mathbb{R}^{n-1}$). If we leave the $a_i$'s undetermined, then $|A+A| = n(n+1)/2$ and $|A-A| = n(n-1)+1$. As $|A-A|$ is larger (in the case where the $a_i$'s are undetermined), in order for $A$ to be an MSTD set we must have non-trivial coincidence of differences, specifically $a_i-a_j=a_k-a_\ell$ for some $(i,j) \neq (k,\ell)$. Given such an equation we can, by projection onto the orthogonal complement of $\mathbb{R}^{n-1}$ of the subspace $(\overrightarrow{e}_i-\overrightarrow{e}_j)-(\overrightarrow{e}_k
-\overrightarrow{e}_\ell)$ spans, represent elements of $A$ by vectors in $\mathbb{R}^{n-1}$. We recompute $|A+A|$ and $|A-A|$. If $|A+A| \le |A-A|$ we pick another non-trivial identification of elements in $A-A$ and repeat the above method with elements of $A$ now represented as vectors in $\mathbb{R}^{n-3}$.  The computation ends with all MSTD sets of size $n$ whose smallest element is 0. With some additional insights that improve the run-time, the program can check $n=8$ fairly quickly; unfortunately $n=9$ is still open (though Hegarty has results for all MSTD sets of size 9 having an additional property).


\section{A positive percentage of sets are MSTD sets}

As for each $n$ studied very few of the $2^n$ subsets of $\{0,1,\dots, n-1\}$ were found to be sum-dominant, it was reasonable to conjecture that in the limit almost no subsets were sum-dominant. While this conjecture is false, the percentage of sum-dominant sets is so small that this error is understandable.

\begin{thm}[Martin - O'Bryant \cite{MO}] As $n\to\infty$, a positive percentage of subsets of $\{0,\dots,n-1\}$ are sum-dominant. \end{thm}

\begin{figure}
\begin{center}
\scalebox{.65764}{\includegraphics{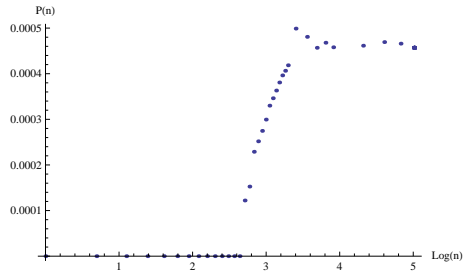}}
\caption{\label{fig:exhaust} The percentage of sum-dominated subsets of $\{0,\dots,n-1\}$ versus $\log n$. These numbers were obtained by enumerating all possible subsets for $n \le 27$, and by simulating 10,000,000 subsets for each $n \in \{30, 35, 40, 45, 50, 75, 100, 125, 150\}$. }
\end{center}\end{figure}

Martin - O'Bryant \cite{MO} proved this probability is at least $2 \cdot 10^{-7}$, which was improved by Zhao \cite{Zh2} to at least $4 \cdot 10^{-4}$; Monte Carlo experiments suggest the true answer is around $4.5 \cdot 10^{-4}$. For small $n$, it is possible to enumerate all subsets of $\{0,\dots,n-1\}$, which we do in Figure \ref{fig:exhaust}.

\begin{figure}
\begin{center}
\scalebox{.5764}{\includegraphics{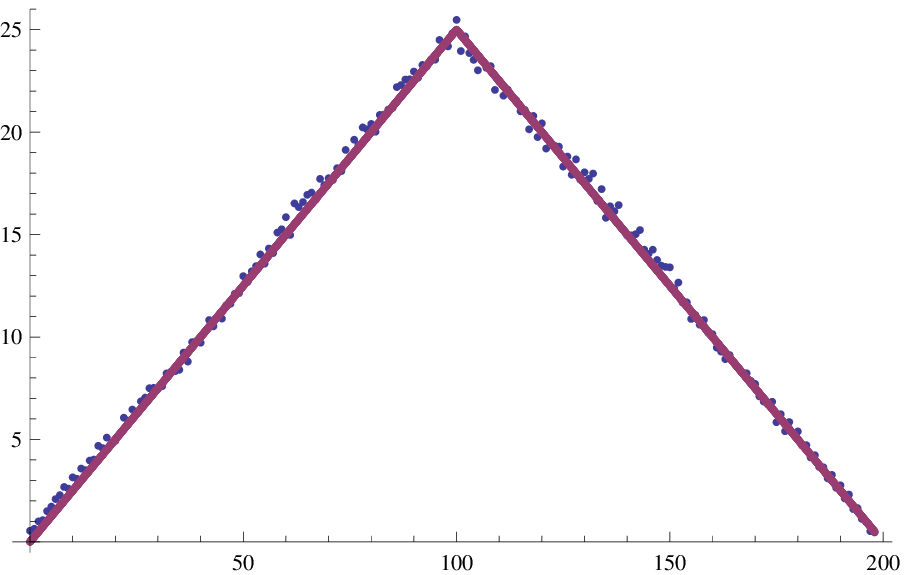}}\ \scalebox{.5764}{\includegraphics{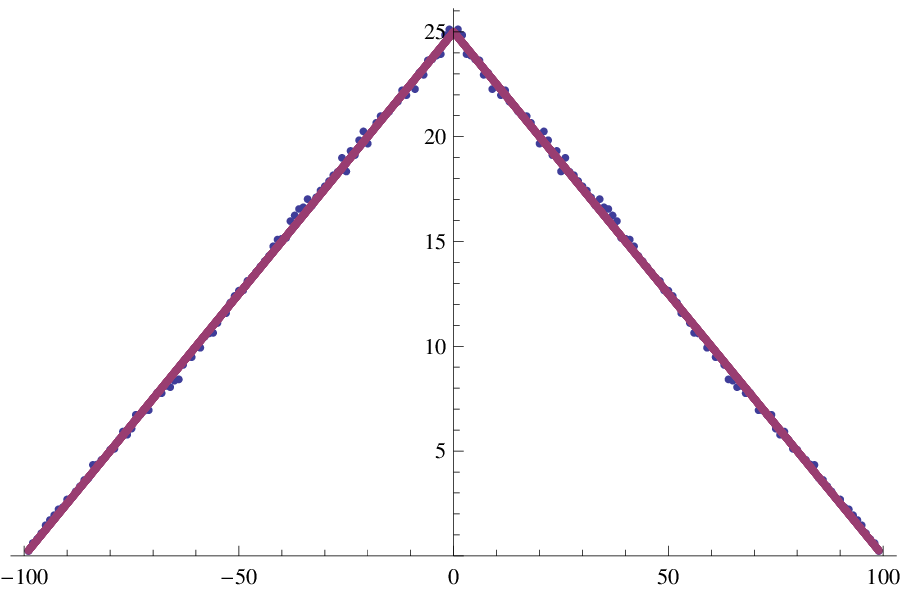}}
\caption{\label{fig:numrepsumdiff} Comparison of predicted and observed number of representations of possible elements of the sumset and difference set for $A \subset \{0,\dots,99\}$ chosen from the uniform model (so each of the $2^{100}$ possible subsets are equally likely to be chosen). We chose 100 different such $A$ and calculated the average number of representations of each possible sum (left plot, which lives in $\{0,\dots,198\}$) and difference (right plot, which lives in $\{-99,\dots,99\}$), compared with the theoretical predictions. Note the spike at 0 was removed from the difference plot. }
\end{center}\end{figure}

Martin and O'Bryant's proof uses probabilistic techniques to estimate the chance that elements are in the sumset and the difference set. For definiteness, consider subsets $S$ of $\{0,1,\dots, n-1\}$. The sumset $S+S$ lies in $\{0,1,\dots, 2n-2\}$ and the difference set $S-S$ in $\{-(n-1),\dots,n-1\}$. The number of representations of a  typical $k \in \{0,1,\dots, 2n-2\}$ as a sum of two elements of $S$  is roughly $n/4 - |n-k|/4$, while the number of representations of a typical $k \in \{-(n-1),\dots,n-1\}$ as a difference of two elements of $S$ is roughly $n/4 - |k|/4$. To see this, first consider the special case when $S = \{0,1,\dots,n\}$. If we want $k = x+y$ with $x\le y$, note once $x$ is chosen then $y$ is determined. If $k \le n-1$ there are essentially $k/2$ choices for $x$; the other case is handled similarly. Our answer differs from $n/4 - |n-k|/4$ by a factor of 2. This factor is due to the fact that a typical set $S$ has approximately $n/2$ elements, and not $n$ elements (by the Central Limit Theorem, the probability is vanishingly small that $|S|$ differs from $n/2$ by more than $n^{1/2+\epsilon}$). Figure \ref{fig:numrepsumdiff} demonstrates the rapidity of convergence. There we uniformly choose many $A \in \{0,\dots,99\}$ and calculate the average number of representations for all the possible sums and differences, and compare with the predictions above. Note for the difference plot we have removed the spike at 0, as for each $A$ there are $|A|$ ways of representing 0 from $A-A$, and by the Central Limit Theorem $|A|$ is approximately $100/2$ or 50.

We see from the above that there are \emph{many} ways to represent the possible sums or differences, so long as they are not near the fringe elements. Their proof proceeds as follows. Let $A$ be an MSTD set, and write $A$ as a disjoint union $L \cup U$, with $L \subset \{0,\dots,\ell-1\}$  and $R\subset \{\ell,\dots, \ell+u-1\}$. Consider the sets $A_M = L \cup M \cup U'$, where $M \subset \{\ell, \dots,\ell + m-1\}$ and $U' = U + m$ (so $U'$ is just $U$ translated by $m$). If $k$ is close to 0 (respectively $\ell+m+u$), then whether or not $k \in A_M + A_M$ depends only on $L + L$ (respectively $U'+U'$). Similarly, the fringe elements of $A_M-A_M$ are determined by $U'-L$ and $L-U'$. By cleverly choosing $A$ (they take $L = \{0,2,3,7,8,9,10\}$ and $U = \{11, 12, 13, 14, 16, 19, 20, 21\}$) we can ensure that there are more sum fringe elements included than difference fringe elements. The proof is completed by showing that  a positive percentage of the possible $M$'s lead to no missing sums or differences in the remaining intervals. This is accomplished through a series of technical lemmas. The estimates here are far from optimal, but suffice to prove a positive percentage of subsets are sum-dominant. Specifically, the authors frequently appeal to the crude estimate that $${\rm Prob}(\{a, a+1, \dots, b\} \nsubseteq A+A) \ \le \ \sum_{k=a}^b {\rm Prob}(k \not\in A+A)$$ (and similarly for difference sets).

There are many other results in this paper. The authors prove the existence of positive lower bounds for the percentage of sum-dominant, balanced, and difference-dominated sets. Though they cannot show the limits exist, they conjecture that this is the case. They show that the average cardinality of the difference sets is four more than the average cardinality of the sumsets, providing additional support that sum-dominant sets should be rare. They also explore $|A+A| - |A-A|$, and show that for any $x$ there is an $A$ such that $|A+A|-|A-A|=x$ with $A \subset \{0,\dots, 17|x|\}$ (which is significantly more economical than the base expansion method would give). The paper ends with some numerical explorations of missing sums, and conjectures that the proportion of subsets $A$ of $\{0,\dots,n-1\}$ with $|A+A| = j$ and $|A-A| = k$ converges to a limiting proportion $\rho_{j,k}$ as $n\to\infty$.

Martin and O'Bryant fixed the fringe (their $L$ and $U$) and varied the middle $M$; Zhao \cite{Zh2} allowed the fringe to vary as well. His methods allow him to obtain MSTD sets that are not missing any middle sums, which he shows happens a vanishingly small number of times. This leads to a significant strengthening of the results of Martin and O'Bryant, and a proof of many of their (and others) conjectures. Specifically, he shows the following limits exist (and provides a deterministic algorithm to approximate their values): the percentage of sets that are sum-dominant; the percentage of sets that are balanced; the percentage of sets that are difference-dominant; the percentage of sets that are missing exactly $s$ sums and $d$ differences; the percentage of sets that have exactly $x$ more sums than differences. The paper ends with an investigation of the probabilities of various elements being in an MSTD set, proving a conjecture of Miller, Orosz and Scheinerman \cite{MOS} that as $n$ grows the probability a `middle' element is in an MSTD set in $\{0,\dots,n\}$  tends to 1/2.


\section{When a `typical' subset is difference-dominated}\label{sec:whentypicalsumdom}

The proofs that a positive percentage of subsets of $\{0,\dots,n-1\}$ are sum-dominant all use, in one way or another, the following fact: if $A$ is uniformly drawn from the $2^n$ subsets of $\{0,\dots,n-1\}$, then with high probability $A$ has essentially $n/2$ elements and almost all possible sums and differences are realized. Along these lines, Martin and O'Bryant \cite{MO} showed that a typical difference set is missing only 7 of the possible differences, and a typical sumset is missing 11 (see \cite{ILMZ} for a proof that the moments of the limiting distribution exist and the tail probabilities are bounded above and below by exponentially decaying probabilities). These techniques apply to a slightly more general case. We may reinterpret the uniform model above as saying each element $k \in \{0,\dots,n-1\}$ is in a subset $A$ with probability $1/2$. We could instead fix a probability $p \in (0,1)$ and let each $k$ be in $A$ with probability $p$.

In this constant probability model, our previous results on a positive percentage again hold. If, however, we allow $p$ to vary with $n$, then the situation is drastically different. Hegarty and Miller \cite{HM1} consider a binomial model where each $k \in \{0,\dots,n-1\}$ is independently chosen to be in a subset $A$ with probability $p(n)$. If $p(n)$ is a constant independent of $n$, we are in the regime handled by Martin and O'Bryant (though we described their method in the uniform model case, similar arguments work so long as the probability is independent of $n$). If, however, $p(n)$ tends to zero, then we are no longer in the case where $|A|$, $|A+A|$ and $|A-A|$ are always large. In this case very few sets are sum-dominant, which is in line with Nathanson's (and others) intuition that, if properly counted, sum-dominant sets are rare.

Before stating their main result, we first set some notation. Let $\mathbb{N}$ denote the positive integers. We say $f(x) = o(g(x))$ if $|f(x)/g(x)| \to 0$ as $x \to \infty$.

\begin{thm}[Hegarty - Miller \cite{HM1}]\label{thm:mainuniform}
Let $p: \mathbb{N} \rightarrow (0,1)$ be any function such that
\be\label{eq:old13} n^{-1}\ =\ o(p(n)) \ \ \ \ {\rm and}\ \ \ \
p(n)\ =\ o(1).\ee For each $n \in \mathbb{N}$ let $A$ be a random
subset of $\{0,\dots, n-1\}$ chosen according to a binomial distribution with
parameter $p(n)$ (so each $k \in \{0,\dots,n-1\}$ is in $A$ with probability $p(n)$). Then, as $n \rightarrow \infty$, the probability
that $A$ is difference-dominated tends to one.
\par More precisely, let $\mathcal{S},
\mathcal{D}$ denote respectively the random variables $|A+A|$ and $|A-A|$.
Then the following three situations arise:
\\
\\
(i) $p(n) = o(n^{-1/2})$ : Then \be\label{eq:old14} \mathcal{S}\ \sim\
{(n\cdot p(n))^{2} \over 2} \;\;\; {\hbox{and}} \;\;\; \mathcal{D}
\sim 2\mathcal{S}\ \sim\ (n \cdot p(n))^{2}. \ee (ii) $p(n) = c \cdot
n^{-1/2}$ for some $c \in (0,\infty)$ : Define the function $g :
(0,\infty) \rightarrow (0,2)$ by \be g(x)\ :=\ 2\left(\frac{e^{-x} -
(1-x)}{x}\right). \ee Then \be\label{eq:old16} \mathcal{S}\ \sim\
g\left({c^{2} \over 2}\right) n \;\;\; {\hbox{and}} \;\;\;
\mathcal{D}\ \sim\ g(c^{2}) n. \ee (iii) $n^{-1/2} = o(p(n))$ : Let
$\mathcal{S}^{c} := (2n+1) - \mathcal{S}$, $\mathcal{D}^{c} :=
(2n+1) - \mathcal{D}$. Then \be \mathcal{S}^{c}\ \sim\ 2 \cdot
\mathcal{D}^{c}\ \sim\ {4 \over p(n)^{2}}. \ee
\end{thm}

The proof proceeds by using various tools to obtain strong concentration results on the sizes of the sum and difference sets. The tools needed depend on the decay of $p(n)$. Not surprisingly, the faster $p(n)$ decays the easier it is to obtain the needed concentration results. The greater the decay, the fewer elements are in a typical $A$, and thus the greater the effect of the non-commutativity of subtraction in generating more new elements. Chebyshev's Theorem suffices for case (i), two still follows elementarily (via a second moment argument), while the third case requires some recent results on strong concentration by Kim and Vu \cite{KiVu,Vu1,Vu2}.

The idea of the proof, at least in case (i), is fairly straightforward. When $n^{-1} = o(p(n))$ and $p(n)= o(n^{-1/2})$, then the expected size of a randomly chosen $A$ is $n p(n) = o(n^{1/2})$. The heart of the proof is to show that such sets are nearly Sidon sets, which means that most pairs of elements generate distinct sums and differences from other pairs (other than the diagonal pairs, those where the two elements are equal, which give just one difference, namely zero). As the non-diagonal pairs generate one sum but two differences, we expect that the difference set will be twice as large as the sumset. A simpler proof of this case is given in the arXiv version of \cite{HM1}, as well as \cite{HM2} (see Appendix 2).

We sketch the proof of case (i) as it highlights the ideas without too many technicalities. The first step is to bound, with high probability, the size of a subset $A$ of $\{0,\dots,n-1\}$ chosen from the binomial model with parameter $p(n) = o(n^{-1/2})$. For ease of exposition, assume $p(n) = c n^{-\delta}$ for some $\delta \in (1/2, 1)$. Using indicator random variables $X_0, \dots, X_{n-1}$ to denote whether or not $k \in A$, by Chebyshev's theorem the probability $X = X_0 + \cdots + X_{n-1}$ is in $[\frac12 cn^{1-\delta}, \frac32 cn^{1-\delta}]$ is at least $1-\frac4{c} N^{\delta-1}$. From here, we obtain upper and lower bounds for the number of pairs of elements $(m,n)$ with $m < n$ both in $A$. All that remains is to show that, with high probability, almost all of the pairs generate distinct sums and differences from each other.

For definiteness we study the differences. If $(m,n)$ and $(m',n')$ generate the same difference then $m-n=m'-n'$. Let $Y_{m,n,m',n'}$ be 1 if $m,n,m',n'$ are in $A$ and $m-n=m'-n'$, and let $Y$ be the sum of the $Y_{m,n,m',n'}$'s. What is $\mathbb{E}[Y]$? Rather than determining it exactly, it suffices to obtain an upper bound. One can show $\mathbb{E}[Y] \le 2C^4 n^{3-4\delta}$ where $C = \max(1,c)$ by considering separately the cases where all four indices are distinct and when three are. As a typical $A$ has size on the order of $n^{1-\delta}$, we expect on the order of $2n^{2-2\delta}$ differences; this is significantly larger than $\mathbb{E}[Y]$, so most of the differences are distinct from each other. All that remains is to control the variance of $Y$, and then another application of Chebyshev's theorem proves that $Y$ is concentrated near its mean, and hence there are on the order of $2n^{2-2\delta}$ differences. The variance estimate follows from elementary counting.

\begin{figure}
\begin{center}
\scalebox{.5764}{\includegraphics{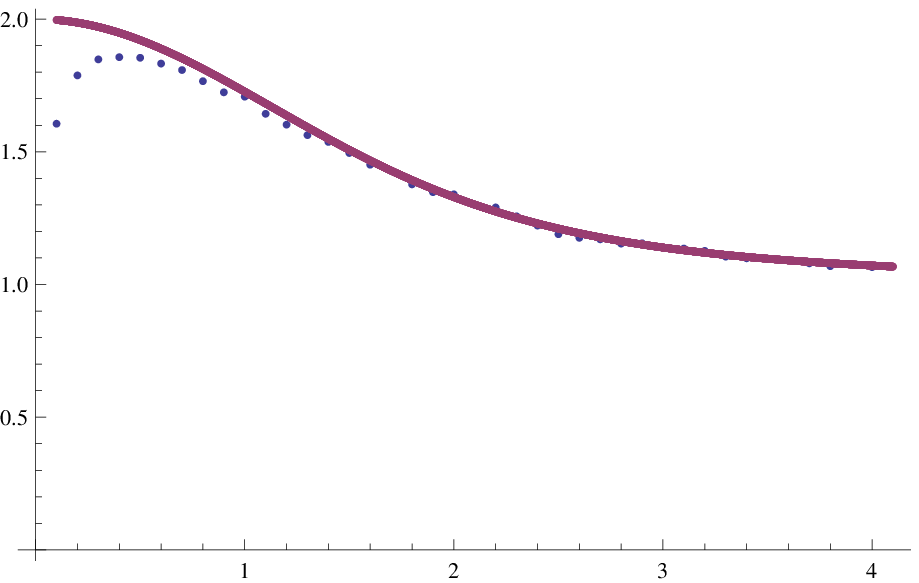}} \scalebox{.5764}{\includegraphics{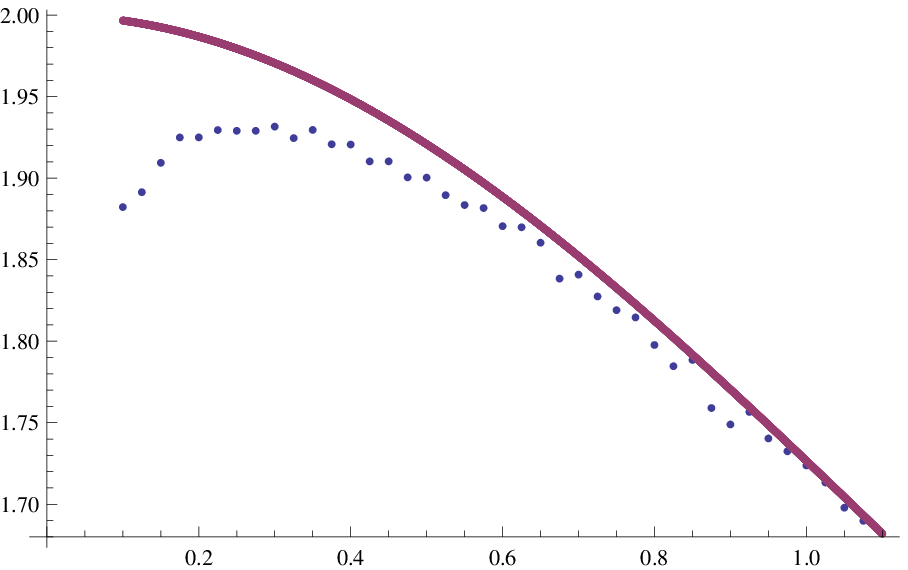}}
\caption{\label{fig:phasetransition1} Plot of $|A-A|/|A+A|$ for ten $A$ chosen uniformly from $\{1,\dots,n\}$ ($n=10,000$ on the left and $100,000$ on the right) with probability $p(n) = c / \sqrt{n}$ versus $g(c^2)/g(c^2/2)$. }
\end{center}\end{figure}

A particularly interesting feature of the above theorem is the existence of a threshold function for the density. If the density $p(n) = o(n^{-1/2})$ then almost surely the ratio of the size of the difference set to the sumset is 2, while above the threshold (so $n^{-1/2} = o(p(n))$) the ratio is 1 (though the number of missing sums is twice that of the number of missing differences). If $p(n) = c n^{-1/2}$ then the ratio of $|A-A|/|A+A|$ tends to $g(c^2)/g(c^2/2)$, with $g(x) = 2\left(e^{-x} - (1-x)\right) / x$. Note this ratio tends to 2 as $c\to 0$ and tends to 1 as $c\to\infty$, which is in line with Cases (i) and (iii) of the theorem. There is thus a nice phase transition in behavior, though this is hard to see experimentally as $10^{-10} n^{-1/2}$ is smaller than $n^{-1/2} \log^{-1} n$ until $n$ exceeds $\exp(10^{10})$. In Figure \ref{fig:phasetransition1} we numerically explore this transition.

\begin{figure}
\begin{center}
\scalebox{.5764}{\includegraphics{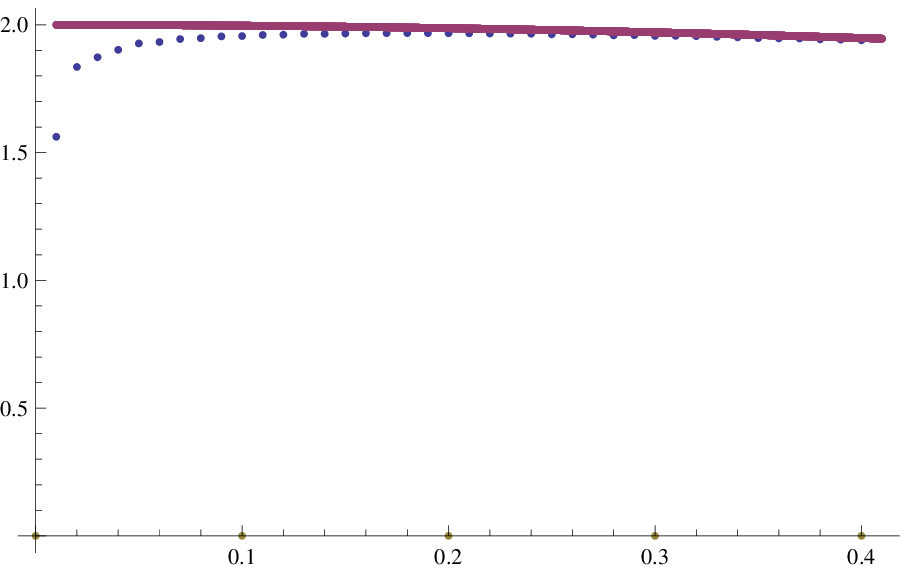}} \scalebox{.5764}{\includegraphics{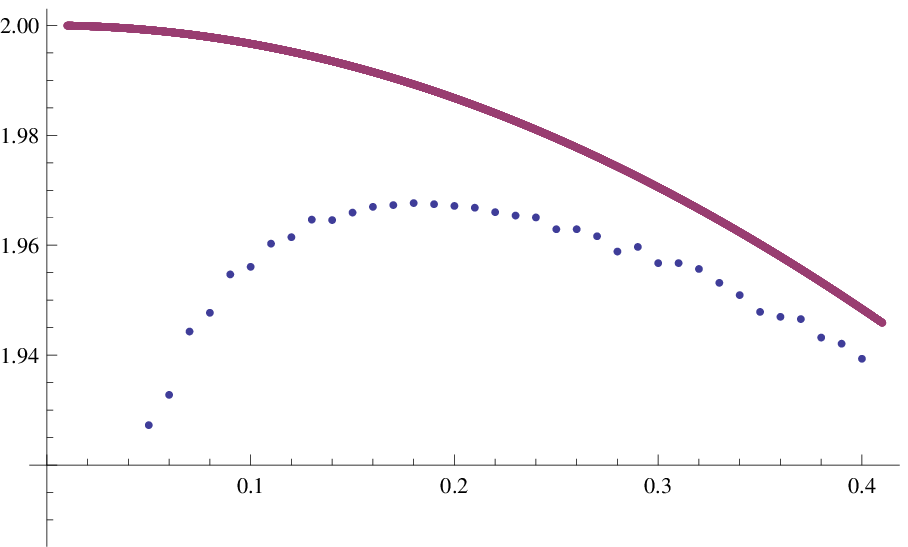}}
\caption{\label{fig:phasetransition2} Plot of $|A-A|/|A+A|$ for ten $A$ chosen uniformly from $\{1, \dots, n\}$ with probability $p(n) = c / \sqrt{n}$ ($n = 1,000,000$) versus $g(c^2)/g(c^2/2)$ (second plot is just a zoom in of the first). }
\end{center}\end{figure}

Not surprisingly, for a fixed $n$ the larger $c$ is, the closer the behavior is to the limiting case. To investigate this further, in Figure \ref{fig:phasetransition2} we examine 40 choices of $c$ from .01 to .41 with $n = 1,000,000$. For $c = .01$ the typical random $A$ has only 10 elements; this increases to about 400 when $c = .41$. We see a noticeable improvement between the observed and conjectured behavior for this larger value of $n$.

To further investigate the transition behavior, we fixed two values of $c$ and studied the ratio for various $n$. We chose $c =.01$ (where the ratio should converge to 1.99997) and $c=.1$ (where the ratio should converge to 1.99667); the results are displayed in Table \ref{ref:tabletwocnvaries}.
\begin{center}
\begin{table}
\begin{center}
\begin{tabular}{|r||c|c|}
  \hline
$n$ & Observed Ratio ($c = .01$) & Observed Ratio ($c = .1$) \\
\hline
    100,000 & 1.123 & 1.873 \\
  1,000,000 & 1.614 & 1.956 \\
 10,000,000 & 1.871 & 1.984 \\
100,000,000 & 1.960 & 1.993\\
  \hline
\end{tabular}
\end{center}
\caption{Observed ratios of $|A-A|/|A+A|$ for $A$ chosen with the binomial model $p(n) = c n^{-1/2}$ for $k \in \{0,\dots,n-1\}$ for $c=.01$ and .1; as $n\to\infty$ the ratios should respectively converge to 1.99997 and 1.99667. Each observed data point is the average from 10 randomly chosen $A$'s, except the last one for $c=.1$ which was for just one randomly chosen $A$.}\label{ref:tabletwocnvaries}
\end{table}
\end{center}


\section{Explicit constructions of large families of MSTD sets}

Until recently, all explicit constructions of families of MSTD sets led to very sparse families, with an exponentially small percentage of the $2^n$ subsets of $\{0,\dots,n-1\}$ being sum-dominant. While the methods of Martin and O'Bryant proved that a positive percentage of the $2^n$ subsets were sum-dominant, their probabilistic method did not allow them to explicitly list these MSTD sets. We quickly review their construction, which was described in greater detail in \S\ref{sec:whentypicalsumdom}.

The word explicit requires some comment. We say a construction is explicit if there is a very simple rule that can quickly be implemented to generate the sets. For example, one method involves taking any set $M \in \{0,\dots, m-1\}$ such that there are never $k$ consecutive elements in $\{0,\dots,m-1\}$ not in $M$. It is very easy to write down sets having this property; it is also easy to count how many such sets there are (and it is this ease in counting that leads to many good results).

Martin and O'Bryant began by choosing a special set $A = L \cup U$ with $L \subset \{0,\dots,\ell-1\}$ and $U \subset \{\ell, \dots, \ell+u-1\}$ such that more of the fringe sums were realized in $A+A$ than fringe differences. They then showed that one could insert almost any set in the middle of $A$ (shifting the elements of $U$ up) and have a sum-dominant set. Miller, Orosz and Scheinerman \cite{MOS} explored which sets, when inserted, did not lead to sum-dominant sets. While this is a very hard question, it turns out that if one carefully chooses sets $L$ and $U$ then one can show \emph{any} set that is never locally too sparse may be inserted and yield a sum-dominant set. The end result is a sparser family than Martin and O'Bryant; however, it is still a large family, and all the technical probability lemmas of \cite{MO} are replaced with elementary counting arguments.

The following property is crucial in the argument. We say a set of integers $A$ has the \emph{property $P_n$} (or is a \emph{$P_n$-set}) if both its sumset and its difference set contain all but the first and last $n$ possible elements (and of course it may or may not contain some of these fringe elements). Explicitly, let $a=\min{A}$ and $b=\max{A}$. Then $A$ is a $P_n$-set if \bea\label{eq:beingPnsetsum} \{2a+n,\ \dots,\ 2b-n\}  \ \subset\  A+A \eea and \bea\label{eq:beingPnsetdiff} \{-(b-a)+n,\ \dots,\ (b-a)-n\}\ \subset\ A-A.\eea It is not hard to show that for fixed $\alpha\in (0,1/2)$ a random set drawn from $\{0,\dots, n-1\}$ in the uniform model is a $P_{\lfloor \alpha n\rfloor}$-set with probability approaching $1$ as $n\to\infty$; it is even easier in our situation as the length of the set $A$ will grow but $n$ will remain fixed. Their main result is

\begin{thm}[Miller-Orosz-Scheinerman \cite{MOS}]\label{thm:mainconstruction}
  Let  $A=L\cup R$ be a $P_n$, MSTD set where  $L\subset \{0,\dots,n-1\}$, $R\subset\{n,2n-1\}$, and $0,2n-1\in A$;\footnote{Requiring $0, 2n-1 \in A$ is quite mild; we do this so that we know the first and last elements of $A$.} for example, $A = \{0, 1, 2, 4, 7, 8, 12, 14, 15\}$ from \cite{Ma} works. Fix a $k \ge n$ and let $m$ be arbitrary. Let $M$ be any subset of $\{n+k,\dots, n+k+m-1\}$ with the property that it does not have a run of more than $k$ missing elements (i.e., for all $\ell \in \{n+k, \dots, n+m\}$ there is a $j \in \{\ell-1,\dots, \ell+k-2\}$ such that $j\in M$). Assume further that $n+k \not\in M$ and set $A(M;k)=L \cup O_1 \cup M \cup O_2 \cup R'$, where $O_1=\{n,\dots, n+k-1\}$, $O_2=\{n+k+m,\dots, n+2k+m-1\}$ (thus the $O_i$'s are just sets of $k$ consecutive integers), and $R'=R+2k+m$. Then

  \ben \item $A(M;k)$ is an MSTD set, and thus we obtain an infinite family of distinct MSTD sets as $M$ varies;

  \item there is a constant $C > 0$ such that as $r\to\infty$ the proportion of subsets of $\{0,\dots,r-1\}$ that are in this family (and thus are MSTD sets) is at least $C / r^4$.

  \een
\end{thm}

It turns out that being a $P_n$-set is not an especially harsh condition, and it is possible to find these sets. The idea of the construction is to add sets in the middle such that all possible middle sums and differences are obtained, and thus whether or not $A(M,k)$ is sum-dominant will depend only on $A$. Specifically, it will depend on whether or not $A$ itself is an MSTD set. While the choices in the construction are not optimal, they do suffice to almost give a positive percentage of sets are sum-dominant, where now we miss by a power instead of by an exponential. A little algebra shows that if $A$ is a $P_n$-set, then so too is our $A(M;k)$. To see this, we need only show that we hit all possible sums and differences except at the fringe. Briefly, the idea behind the construction is that because $O_1$ and $O_2$ have $k$ consecutive integers and $M$ never misses $k$ consecutive integers, when we look at sums such as $O_1 + M$ we will always have two elements in $A(M;k)$ that will add to the desired number (and similarly for the differences).

The rest of the proof deals with examining how restrictive the assumption is that $M$ never misses $k$ consecutive integers. One can solve this by writing down a recurrence relation, but an elementary approach is available which yields quite good results with little work. We assume a slightly stronger condition: we break $M$ into blocks of length $k/2$ and assume $M$ always has an element from each of these blocks. This ensures that there can never be a gap as large as $k$ between elements of $M$ (the gap is at most $k-2$). There are $2^{k/2}$ possibilities for each block of length $k/2$; all but one (choosing no elements) satisfies the stronger condition. The percentage of such valid sets in $\{0,\dots, r-1\}$ is a constant times \be\label{eq:keycardsum} \sum_{k=n}^{r/4} \frac1{2^{2k}} \left(1 - \frac1{2^{k/2}}\right)^{\frac{r}{k/2}}. \ee There are two factors leading to obtaining less than a positive percentage. The first is, obviously, that in each block of length $k/2$ we lose one possibility, and this factor is raised to a high power. The second is that $O_1$ and $O_2$ are completely determined and their length depends on $k$. Thus, as soon as $k$ grows with $n$, we see we cannot have a positive percentage. Analyzing the sum gives the claimed bounds.

\begin{rek}\label{rek:schilling} The above theorem can be improved by appealing to an analysis of the probability $m$ consecutive tosses of a fair coin has its longest streak of consecutive heads of length $\ell$ (see \cite{Sc}). What is fascinating about the answer is that while the expected value of $\ell$ grows like $\log_2(m/2)$, the variance converges to a quantity independent of $m$, implying an incredibly tight concentration. If we take $O_1$ and $O_2$ as before and of length $k$, we may take a positive percentage of all $M$'s of length $m$ to insert in the middle, so long as $k = \log_2(m/2) - c$ for some $c$. The size of $A$ is negligible; the set has length essentially $m+2k$. Of the $2^{m+2k}$ possible middles to insert, there are $C 2^m$ possibilities (we have a positive percentage of $M$ work, but the two $O$'s are completely forced upon us). This gives a percentage on the order of $2^m / 2^{m+2k}$; as $k=\log_2(m/2)-c$, this gives on the order of $1/m^2$ as a lower bound for the percentage of sum-dominated sets, much better than the previous $1/m^4$.
\end{rek}

The results of \cite{MOS} can be generalized to compare linear forms. We can find infinite families of sets satisfying \be\label{eq:genMOS} \left|\gep_1 A + \cdots + \gep_n A\right| \ > \ \left|\widetilde{\gep}_1 A + \cdots + \widetilde{\gep}_n A\right|, \ \ \ \gep_i, \widetilde{\gep}_i \in \{-1,1\} \ee \emph{if} we can find one set satisfying the above. We've seen from \cite{MO,Zh2} that very few sets are sum-dominant; thus we expect the percentage of sets satisfying \eqref{eq:genMOS} to be extremely small, and thus expect it to be a challenge to find the needed set. Brute force search found $\{$0, 1, 2, 3, 7, 11, 17, 21, 22, 24, 25, 28, 29, 30, 31, 33,    44, 45, 48, 49$\}$, which gives $|A+A+A| > |A+A-A|$; unfortunately, such naive searching was unsuccessful in finding examples for other comparisons. We describe a new method by Iyer, Lazarev, Miller and Zhang \cite{ILMZ} in \S\ref{sec:generalizedMSTDsets} which generates the needed sets to begin the induction arguments.

In the above generalizations, the construction from \cite{MOS} with $|A+A| > |A-A|$ is mimicked for the linear forms. In particular, we still assume that $M$ has at least one element in each block of length $k/2$. While this was necessary for $|A+A| > |A-A|$, Miller, Pegado and Robinson \cite{MPR} show that this is not needed in general. For example, if we are studying $|A+A+A+A|$ versus $|A+A-A-A|$, we are assisted by the fact that we can have $O_i + O_j$ and then add this to $M+M$. The final result of all of this is that we may allow $O_1$ and $O_2$ to be significantly more sparse than in \cite{MOS}, where they had to choose $k$ consecutive elements and thus had no freedom. What matters is that $O_i+O_j$ contain large consecutive blocks of integers, not that each $O_i$ do so. This allows us to improve upon the $1/2^{2k}$ terms in \eqref{eq:keycardsum}.

Before stating the result, we need to slightly generalize the notion of a $P_n$-set to a $P_n^4$-set. We say $A$ is a $P_n^4$-set if $A + A + A + A$ and $A + A - A - A$ each contain all but the first and last $n$ elements; thus what we called a $P_n$-set before is really a $P_n^2$-set.\\

\begin{thm}[Miller-Pegado-Robinson \cite{MPR}] Let  $A=L\cup R$ be a $P_n$, MSTD set where  $L\subset \{0,\dots,n-1\}$, $R\subset\{n,2n-1\}$, and $0,2n-1\in A$;\footnote{As before, requiring $0, 2n-1 \in A$ is quite mild and is done so that we know the first and last elements of $A$.} for example, $A =  \{$0, 1, 3, 4, 7, 26, 29, 30, 32, 33, 34, 27, 28, 31, 53, 56, 57, 59, 60, 61$\}$ works. Fix a $k \ge n$ and let $m$ be arbitrary. Let $M$ be any subset of $\{n+k,\dots, n+k+m-1\}$ with the property that it does not have a run of more than $k$ missing elements (i.e., for all $\ell \in \{n+k, \dots, n+m\}$ there is a $j \in \{\ell-1,\dots, \ell+k-2\}$ such that $j\in M$). Assume further that $n+k \not\in M$ and set $A(M;k)=L \cup O_1 \cup M \cup O_2 \cup R'$, where $O_1=\{n,\dots, n+k-1\}$, $O_2=\{n+k+m,\dots, n+2k+m-1\}$ (thus the $O_i$'s are just sets of $k$ consecutive integers), and $R'=R+2k+m$. Then

  \ben \item $A(M;k)$ is an MSTD set, and thus we obtain an infinite family of distinct MSTD sets as $M$ varies.

  \item There is a constant $C > 0$ such that as $r\to\infty$ the proportion of subsets of $\{0,\dots,r-1\}$ that are in this family (and thus are MSTD sets) is at least $C / r^{4/3}$.

  \item With better choices of $O_1$ and $O_2$, one can explicitly construct a large family of sets $A$ with $|A+A+A+A| >
|(A+A)-(A+A)|$ and show that the density of sets $A \subset \{0,\dots, n-1\}$
satisfying this condition is at least $C/n^r$, where $r = \frac16
\log_2(256/255) \le .001$.

\item For each integer $k$, there is a set $A \subset \{0, \dots, 157k\}$ such that $|2A+2A|$ $-$ $|2A-2A|$ $= k$; if $k$ is large we may take $A \subset \{0, \dots, 35|k|\}$.

  \een
\end{thm}


The proof of the first two assertions follows identically as in \cite{MOS} (if we argue as in Remark \ref{rek:schilling} and use the results from \cite{Sc}, we may improve (2) from $r^{4/3}$ to $r^{2/3}$). For the third assertion, the additional binary operations gives us enormous savings and removes many of the restrictions on the form of the $O_i$'s. We note that the $O_i$'s show up in sums and differences at least in pairs, unless matched with $L+L+L$, $R'+R'+R'$ or $L+L-R'$ ($A = L\cup R$). Each of $L+L+L$, $R'+R'+R'$ and $L+L-R'$ contains a run of 16 elements in a row for our set $A$. This allows us to relax the restrictions on $O_i$ from \cite{MOS} (each $O_i$ was $k$ consecutive elements); if each $O_i$ has no run of 16 missing elements
and $2O_i$ is full for both $O_i$'s, simple algebra shows that we get all sums and differences as before. This looser structure on the $O_i$'s  allows us to replace the $1/2^{2k}$ in \eqref{eq:keycardsum} with a much better term, leading to a significantly better exponent and thus greatly improve the density bound.


Returning to MSTD sets (and not their generalizations), the current record for densest explicit family of MSTD sets is due to Zhao \cite{Zh1}, who found a family of $\{0, \dots, n-1\}$ of order $2^n/n$. He achieved this by showing a correspondence between bidirectional ballot sequences and sum-dominant sets. A \emph{ballot sequence} is a list of 1s and 0s if every prefix has more 1s than 0s and the maximum excess of 1s over 0s is attained at the end of the sequence. If you imagine the 1s as winning \$1 and the 0s as losing \$1, we may interpret this as we bet a fixed amount each game, our winnings are always positive and our greatest balance is at the end. A sequence of 1s and 0s is a \emph{bidirectional ballot sequence} if both it and the reversed sequence are ballot sequences.

Much of the construction is similar to \cite{MO,MOS}; we again take a set that leads to the desired fringe behavior, and study which sets $M$ may be inserted. Unlike the previous constructions, here we ask that $M$ is a bidirectional ballot sequence (where we write 1 if an element is in $M$ and 0 if it is not). This is equivalent to the following. Let $M \subset \{0,\dots,m-1\}$. Then every prefix and suffix of $\{0,\dots,m-1\}$ has more than half its elements in $M$. As each prefix and suffix has more than half its elements in $M$, by the pidgeon hole principle at least one pair will be in $M$, and that will generate the desired sum or difference. The problem is thus reduced to counting the number of bidirectional ballot sequences,



\section{Generalized MSTD Sets}\label{sec:generalizedMSTDsets}

There are many ways to generalize the notion of a sum-dominant set. Below we discuss two possibilities that were recently analyzed in \cite{ILMZ}; we comment briefly on the ideas and constructions, and refer the reader to the article for full details. As we are always adding sets and never multiplying, in all arguments below we use the shorthand notation \be
kA \ = \ \underbrace{A+\cdots+A}_{k\text{ times }}. \ee 

\begin{enumerate}

\item Given non-negative integers $s_1, d_1, s_2, d_2$ with $s_1+d_1 = s_2+d_2 \ge 2$, can we find a set $A$ with $|s_1A-d_1A| > |s_2A-d_2A|$, and if so, does this occur a positive percentage of the time?
    
\item We say a set is $k$-generational if $A$, $A+A$, $\dots$, $kA$ are all sum-dominant. Do $k$-generational sets exist, and if so, do they occur a positive percentage of the time? Is there a set that is $k$-generational for all $k$?
    
\end{enumerate}

The first question is motivated by generalizing the binary comparison. When $s_1+d_1 = 2$, the only possible sets are $A+A$ and $A-A$ (note $-A-A$ is the same as the negation of $A+A$). When $s_1+d_1=3$, again there are again essentially just two possibilities, $A+A+A$ and $A+A-A$ (as $A-A-A = -(A+A-A)$, and thus without loss of generality we might as well assume $s_i \ge d_i$). The situation is markedly different once the sum is at least 4. In that case, we now have $A+A+A+A$, $A+A+A-A$ and $A+A-A-A$. All possible orderings happen a positive percentage of the time.

\begin{thm}[Iyer-Lazarev-Miller-Zhang \cite{ILMZ}]\label{thm:generalizedcomparisons} Given non-negative integers $s_1, d_1, s_2, d_2$ with $s_1+d_1 = s_2+d_2 = k \ge 2$, if $\{s_1,d_1\} \neq \{s_2,d_2\}$ then a positive percentage of all sets $A$ satisfy $\left|s_1A-d_1A\right|>\left|s_2A-d_2A\right|$. For definiteness assume $s_1$ is the largest of the $s$'s and $d$'s. Given any non-negative integers $i,j$ with $j \le 2i$, for all $n$ sufficiently large there exists an $A \subseteq \{0, 1, \dots, n\}$ such that
$|s_1A - d_1A | = kn+1 -i$ and $|s_2A - d_2A| = kn +1 - j$.
\end{thm}

\emph{Sketch of the proof.} The proof is similar in spirit to many of the results in the field; we first find one example by cleverly constructing a set with a certain fringe structure, and then use the methods from Martin-O'Bryant \cite{MO} to expand the set by essentially adding anything in the middle. The difficulty, as was apparent in \cite{MOS}, is in constructing one such set. To make such a set $A$, we pick fringes $L$ and $R$ such that their
sums (with themselves or with each other) have the same structure
(a few chosen elements below the maximum missing). Then we let $A=L\cup M\cup(n-R)$,
where $M$ is a large interval in the middle. If $M$ is large enough,
we don't have to worry about anything besides the fringes. As $A$
is summed, the fringes slowly fill in, however, we choose $L$ such
that $\max(L)<\max(R)$. This means that the right fringe of $kA$
fills in faster than the left. Note that the right fringe of $kA$
is just $k(n-R)$, and the right fringe of $s_{2}A-d_{2}A$ is $s_{2}(n-R)-d_{2}L$.
Since $R$ grows faster than $L$, we can choose the middle such that
$k(n-R)$ will intersect with the middle and be filled in, but $s_{2}(n-R)-d_{2}L$
will not. At the same time, we have that the left fringe of $kA$
is missing one element, and the left fringe of $s_{2}A-d_{2}A$ is
as well. We refer the reader to \cite{ILMZ} for details of the construction for a given $i$ and $j$.

To illustrate the method, consider
\begin{align}
L &\ = \ \{0,1,3,4,\ldots,k-1,k,k+1,2k+1\}\nonumber\\
& \ = \ [0,\ell]\backslash\left(\{2\}\cup[\ell-k+1,\ell-1]\right)\nonumber\\
R & \ = \  \{0,1,2,4,5,\ldots,k,k+1,k+2,2k+2\}\nonumber\\
 &\ =\ [0,r]\backslash\left(\{3\}\cup[k+3,2k+1]\right).\end{align}
For any $x,y\in\mathbb{N}$, the basic structure of $xL+yR$ is the same as that of the original
set. Basically, $xL+yR$ is always missing the first $k$ elements
below the maximum, as well as the singleton element $2k-1$ away from
the maximum. Even more, it is missing no other elements.

Returning to the original problem, our initial set has a fringe structure and sufficient
empty space to allow the fringe to grow and exhibit the desired behavior,
followed by a full middle. We can have more control of the set's behavior
by putting in another fringe along the outside, with sufficient empty
space to let the fringe exhibit the correct behavior before it intersects
with the inner fringe. This process becomes technical, but it allows for a great degree of control over sets. 

More generally, one has

\begin{thm}[Iyer-Lazarev-Miller-Zhang \cite{ILMZ}]\label{thm:kgen} Given finite sequences of length $k$ called $x_{j},y_{j},w_{j},z_{j}$
such that $x_{j}+y_{j}=w_{j}+z_{j}=j$, $x_{j}\neq w_{j}$  and $x_{j}\neq z_{j}$,
for every $2\leq j\leq k$, there exists a set $A$ such that $\left|x_{j}A-y_{j}A\right|>\left|w_{j}A-z_{j}A\right|$
for every $2\leq j\leq k$. In particular, there exists a set $A$ such that $\left|cA+cA\right|>\left|cA-cA\right|$
for every $1\leq c\leq k$.\end{thm}

The above theorem answers our second question, and is the best possible (at least in regard to $k$-generational sets) as every set is finite generational. In other words, one cannot have a set $A$ such that $\left|cA+cA\right|>\left|cA-cA\right|$
for all $c$. It turns out that all sets have
a kind of limiting behavior. As we continue adding $A$ to its sums, eventually we have a full
middle, and any interesting behavior will occur on the fringes. Note that if we normalize $A$ to include $0$, we have $cA\subset cA-cA$.
Essentially, the difference sets eventually have each fringe
element as the sum sets. When $c$ is sufficiently large, the fringes
of $cA$ stabilize, which gives
$\left|cA-cA\right|\geq\left|cA+cA\right|$. Now, taking differences
allows the left fringe to interact with the right fringe, while taking
only sums keeps these separate. This means that it is possible (and
in fact likely) to have $\left|cA-cA\right|>\left|cA+cA\right|$ for
all sufficiently large $c$. We can readily obtain an upper bound on how long we must wait for the limiting behavior of $|kA|$ to set in.

\begin{thm}[Iyer-Lazarev-Miller-Zhang \cite{ILMZ}] Let $A=\{a_1,a_2,\dots,a_m\}\subset$ $\{0, 1, \dots, n-1\}$ be a set of integers ($a_1<a_2<\ldots<a_m$) and let $s=\mbox{gcd}(a_1$, $a_2$, $\dots$, $a_m)$. Then there exists an integer $N$ such that for $k\ge N$ we have $|kA|=\frac{k(a_m-a_1)}{s}-C$ where $C$ is a constant and $k$ is bounded above by $\frac{a_m-a_1}{s}$. \end{thm}

\emph{Sketch of the proof}: It is enough to show the claim for a set of the form $\{0, a_1,\ldots,a_m\}$ with $\mbox{gcd}(a_1,\ldots,a_m)=1$. Adding $A$ to itself $a_1$ times will generate all congruence classes of $a_1$ because of $\mbox{gcd}(a_1,\ldots,a_m)=1$. Adding $A$ to itself $a_m$ times will make both the left ($L$) and right ($R$) fringes stabilize, where $L=kA\cap \{0$, $1$, $\dots$, $a_1a_m\}$ and  $R=kA\cap \{ka_m-a_1a_m, \dots, ka_m\}$, and also ensures that the middle part is completely filled.

\ \\

We end with a few examples of the previous theorems. In these theorems no effort was made to optimize the arguments and generate minimal such sets; this would be an interesting future project, as it is almost surely possible to construct examples of sets with the above
properties that contain many fewer elements. In particular, the base
expansion method of combining sets is extremely inefficient. An alternative,
which is discussed briefly above, is the multiple fringes method.
This allows for much smaller sets, however, the requirements for the
method to work are very stringent, and the proofs are messy. Therefore
we find it best to give the constructions using the base expansion
method instead.

\begin{itemize}
\item If we set \be A \ = \ \{0,1,3,4,5,9,33,34,35,50,54,55,56,58,59,60\} \ee then \be \left|A+A+A+A\right|\ > \ \left|A+A+A-A\right|.\ee

\item If we take \bea A & \ = \ & \{0,1,3,4,7,26,27,29,30,33,37,38,40,41,42,43,46,49,50,52,53,54, \nonumber\\ & & \ \ \ 72,75,76,78,79,80\}\eea then \be \left|A+A\right| \ > \ \left|A-A\right| \ \ \ {\rm and} \ \ \  \left|A+A+A+A\right| \ > \ \left|A+A-A-A\right|; \ee in other words, $A$ is 2-generational.

\item If we let\begin{align}
A\ = \ & \{0,1,3,4,5,6,11,50,51,53,54,55,56,61,97,132,137,138,140,\nonumber\\
 & 142,143,144,182,187,188,189,190,192,193,194\}\end{align} then \be \left|4A-A\right| \ > \ \left|5A\right| \ \ \ {\rm and} \ \ \ \left|4A-A\right|>\left|3A-2A\right|. \ee
 
\end{itemize}


%
%
%

\ \\

\end{document}